\documentclass[journal ]{new-aiaa}
\usepackage[utf8]{inputenc}
\usepackage{textcomp}

\usepackage{graphicx}
\usepackage{amsmath}
\usepackage[version=4]{mhchem}
\usepackage{siunitx}
\usepackage{longtable,tabularx}
\usepackage{tikz}
\usepackage{pgfplots}
\usepackage{pgf}
\usepgfplotslibrary{groupplots}
\pgfplotsset{width=7.5in, compat=1.18}

\usepackage{bm}
\usepackage{bbm}
\usepackage{amsfonts}

\usepackage{amssymb}
\usepackage{amsthm}

\usepackage{algorithm}
\usepackage{algpseudocode}

\usepackage{booktabs}

\usepackage{multirow} 
\usepackage{subcaption}

\usepackage{trimclip}
\def\AR{\clipbox{0pt 0pt .32em 0pt}\AE\kern-.30emR}

\newcommand{\mycomment}[1]{}

\title{Multi-fidelity approaches for general constrained Bayesian optimization with application to aircraft design}

\author{O. Cordelier~\footnote{MSc student, Department of Mathematics and Industrial Engineering, oihan.cordelier@polymtl.ca, AIAA Member.} and Y. Diouane~\footnote{Associate Professor, Department of Mathematics and Industrial Engineering, youssef.diouane@polymtl.ca, AIAA MDO TC Member.}}
\affil{GERAD and Department of Mathematics and Industrial Engineering, Polytechnique Montréal.}

\author{N. Bartoli~\footnote{Senior researcher, DTIS, nathalie.bartoli@onera.fr, AIAA MDO TC Senior Member.}}
\affil{DTIS, ONERA, Université de Toulouse, 31000, Toulouse, France}
\affil{F\'ed\'eration ENAC ISAE-SUPAERO ONERA, Université de Toulouse, France}

\author{E. Laurendeau~\footnote{Full Professor, Department of Mechanical Engineering, Polytechnique Montréal.}}
\affil{Department of Mechanical Engineering, Polytechnique Montréal.}

\begin{document}

\maketitle

\begin{abstract}
Aircraft design relies heavily on solving challenging and computationally expensive Multidisciplinary Design Optimization problems. In this context, there has been growing interest in multi-fidelity models for Bayesian optimization to improve the MDO process by balancing computational cost and accuracy through the combination of high- and low-fidelity simulation models, enabling efficient exploration of the design process at a minimal computational effort. In the existing literature, fidelity selection focuses only on the objective function to decide how to integrate multiple fidelity levels, balancing precision and computational cost using variance reduction criteria. In this work, we propose novel multi-fidelity selection strategies. Specifically, we demonstrate how incorporating information from both the objective and the constraints can further reduce computational costs without compromising the optimality of the solution. We validate the proposed multi-fidelity optimization strategy by applying it to four analytical test cases, showcasing its effectiveness. The proposed method is used to efficiently solve a challenging aircraft wing aero-structural design problem. The proposed setting uses a linear vortex lattice method and a finite element method for the aerodynamic and structural analysis respectively. We show that employing our proposed multi-fidelity approach leads to $86\%$ to $200\%$ more constraint compliant solutions given a limited budget compared to the state-of-the-art approach.
\end{abstract}

\section*{Nomenclature}


{\renewcommand\arraystretch{1.0}
\noindent\begin{longtable*}{@{}l @{\quad=\quad} l@{}}
BO & Bayesian Optimization\\
$C_D$ & Wing Drag Coefficient\\
$C_L$ & Wing Lift Coefficient\\
DoE & Design of Experiment\\
EGO & Efficient Global Optimization\\
GP & Gaussian Process\\
HF & High-fidelity\\
LF & Low-fidelity\\
LHS & Latin Hypercube Sampling\\
MDA & Multidisciplinary Design Analysis\\
MDO & Multidisciplinary Design Optimization\\
MF & Multi-fidelity\\
PoF & Probability-of-Feasibility\\
VF & Variable-fidelity\\
VLM & Vortex Lattice Method\\
\end{longtable*}}

\section{Introduction}\label{sec:1}
\lettrine{M}ultidisciplinary Design Optimization (MDO)~\cite{RaymerAircraftDesignConceptual2018} aims to find the best design by considering the trade-offs and dependencies between  disciplines (e.g., aerodynamics, structural mechanics, and thermodynamics). Multidisciplinary Design Analysis (MDA) is the study of interacting disciplines, where a design variable of a discipline might influence another. MDO builds upon MDA by incorporating it into an optimization process. Multi-fidelity models also play a crucial role in balancing the accuracy and computational efficiency of MDO across the disciplines. High-fidelity models, such as Reynolds-Averaged Navier-Stokes (RANS), are used to accurately simulate complex flow phenomena, but their computational cost makes them suitable for the final stages of MDO when precision is critical. Euler equations, as medium-fidelity models, simplify certain aspects, such as ignoring viscosity, to reduce computational effort while maintaining reasonable accuracy, which is ideal for intermediate stages (e.g., induced drag).  Low-fidelity models, such as those based on Laplace equation (e.g., incompressible flows) provide quick approximations based on simplified assumptions, allowing for rapid exploration of design options across disciplines. By strategically using these models within MDO, designers can efficiently explore and optimize complex designs while ensuring that computational efforts are used effectively. 

MDO problems can be formulated as a computationally expensive-to-evaluate blackbox constrained optimization of the form
\begin{equation}
\begin{aligned}
\min_{\mathbf{x} \in \mathbb{R}^d} \quad & f(\mathbf{x}) \\
\text{s.t.} \quad 
& \bm g(\mathbf{x}) \le \bm 0\\
& \bm h(\mathbf{x}) = \bm 0\\
\end{aligned}
\label{eq:opt_prob}
\end{equation}
where $f:\mathbb{R}^d \mapsto \mathbb{R}$ is the objective function,  \hbox{$\bm g:\mathbb{R}^d \mapsto \mathbb{R}^m$} gives the $m$ inequality constraints, and \hbox{$\bm h:\mathbb{R}^d \mapsto \mathbb{R}^p$} gives the $p$ equality constraints.
The design space $\Omega \subset \mathbb{R}^d$ is a bounded domain. The functions $f$, $\bm g$ and $\bm h$ are typically simulations with no exploitable properties such as structure or derivatives (i.e., blackbox). Calling $f$, $\bm g$, and $\bm h$ is often expected to be very computationally expensive and thus takes a long time to evaluate. In this paper, we assume access to multiple sources of fidelity data related to the objective function $f$ and the constraint functions $\bm g$ and $\bm h$. The fidelity levels are assumed to range from the lowest fidelity (cheapest) to the highest fidelity (most expensive). Examples of challenging blackbox MDO problems are commonly encountered in the aerospace industry. See, for instance,~\citet{priem_2020}, where Bombardier presents an aircraft optimization use case to identify configurations from early-stage concepts to detailed aircraft designs. 

\textit{Bayesian Optimization} (BO)~\cite{frazier2018tutorial,SEGO-UTB} is a powerful strategy for solving expensive blackbox problems as in Eq.~\eqref{eq:opt_prob}. The use of multiple fidelities to improve BO was explored both for mono-objective \cite{meliani_multi-fidelity_2019} and  multi-objective~\cite{charayron2023towards} optimization. These latest BO approaches are largely based on powerful multi-fidelity modeling tools of~\citet{kennedy_predicting_2000}. Building on this foundation,~\citet{ForresterAlexanderI.J2007Movs} proposed a multi-fidelity extension for an unconstrained BO strategy, known as the Efficient Global Optimization (EGO) algorithm. Later, EGO was extended to the constrained setting with the development of the \textit{Super EGO} framework (SEGO) \cite{SasenaExplorationmetamodelingsampling2002}. \citet{meliani_multi-fidelity_2019} proposed MFSEGO, an extension of SEGO to the multi-fidelity setting. In~\cite{meliani_multi-fidelity_2019,charayron2023towards}, a fidelity selection criterion focusing solely on the objective function was introduced to decide how to integrate multiple fidelity levels. The fidelity models related to the constraints did not influence the multi-fidelity optimization strategy.
Other multi-fidelity optimization frameworks exist, including MFSKO~\cite{Huang2006} and NM2-BO~\cite{DiFiore2024}; however, they provide limited discussion on constraint-handling strategies. MFSKO proposes possible ideas to integrate constraints but does not provide benchmark results. NM2-BO addresses constraints by reformulating them as penalty terms added to the objective function. Frameworks such as VF-EI~\cite{Zhang2018} and VF-PI~\cite{Ruan2020} explicitly integrate multi-fidelity constraints, but they cannot efficiently treat equality constraints as they penalize their respective acquisition function using the probability-of-feasibility~\cite{Schonlau1997}.

In this paper, we investigate the potential of incorporating multiple fidelity levels for constraints, in addition to those related to the objective function, into the optimization process. These range from the most optimistic approach, where decisions are based mainly on the best affordable fidelity level, to the most pessimistic approach, where decisions are guided by the worst affordable fidelity level. 
The proposed approaches are then used to solve a bi-fidelity aero-structural wing optimization case.
The obtained results suggest that taking into account both the objective and constraint functions when deciding the level of fidelity to be used within the optimization process is beneficial; the pessimistic strategy appears to perform best in our tests.

The outline of the paper is as follows. 
In Section~\ref{sec:2}, we provide a detailed review of multi-fidelity Gaussian processes. Section~\ref{sec:3} reviews BO and its multi-fidelity extension, MFSEGO. Our proposed methodology is introduced in Section~\ref{sec:4}, where we present a detailed analysis of the properties of the proposed algorithm.  
Section~\ref{sec:4} discusses implementation details and presents the results obtained from four challenging analytical test cases. The results are compared against two other multi-fidelity frameworks, VF-EI and VF-PI. Results from a wing aerostructural design optimization case are also provided.  
Finally, conclusions and ongoing work are summarized in Section~\ref{sec:5}.  

\section{Multi-fidelity Gaussian processes}\label{sec:2}
\subsection{Gaussian processes}
Gaussian Processes (GP)~\cite{rasmussen_gaussian_2008} build surrogate models capable of predicting function $y$ values and the corresponding uncertainty of the prediction across its domain from a limited number of function samples. Using~\citet{sacks_designs_1989} formulation, the function is modeled following a linear interpolation and a stochastic process, which captures the departure from said linear model as written: 

\begin{equation}
y(\bm{x}) = \sum_{k=1}^{q}{\beta_k f_k(\bm{x})} + Z(\bm{x})
\label{eq:gp}
\end{equation}
where \(f_k(\bm{x})\) are $q$ basis functions, and \(\beta_k\) are their corresponding weights. $Z(\bm x)$ is characterized by a covariance function: 
\begin{equation}
\label{eq:correlation_matrix}
\text{cov}\left[ Z(\bm{x}^i), \; Z(\bm{x}^j)\right ] = \sigma^2_z R\left( Z(\bm{x}^i), Z(\bm{x}^j) \right)= \sigma^2_z\prod_{l=1}^{d} \exp \left( -\theta_l \left| x_l^i - x_l^j \right|^2 \right),
\end{equation}
where $R$ is the squared exponential correlation function for $\bm{x}^i$ and $\bm{x}^j$ two points in $\mathbb{R}^d$. 
The scale factor $ \sigma_z $ and the vector of $ d $-dimensions \(\bm{\theta}\) correspond to the prior variance and the correlation lengths $ \theta_l $ in each dimension~\cite{rasmussen_gaussian_2008}.

The Design of Experiments (DoE), denoted as \(\text{DoE} = \{ (\bm{x}, y)_k \}_{k=1, \dots, N} \), gathers \(N\) samples from the function \(y\) from which it is then possible to approximate the hyperparameters by using the Maximum Likelihood Estimation (MLE)~\cite{bachoc2013cross} and  compute vectors and matrices at a query point $\bm{x}$: 
\begin{equation}
\label{eq:gp_structure}
\begin{split}
\bm{f}(\bm{x}) & = (f_1(\bm x),\dots,f_q(\bm x))^\top \in \mathbb R^q\\
\bm{r}(\bm{x}) & = \left( R(\bm{x}_1, \bm{x}), \dots, R(\bm{x}_N, \bm{x}) \right)^\top \in \mathbb{R}^{N}\\
\bm{\beta} & = (\beta_1,\dots,\beta_q)^\top \in \mathbb R^q\\
\bm{X} & = \left( \bm{x}_1, \dots, \bm{x}_N \right)^\top \in \mathbb{R}^{N\times d}\\
\bm{Y} & = \left( y_1, \dots, y_N \right)^\top \in \mathbb{R}^{N}\\
\bm{F} & = \left( \bm{f}(\bm{x}_1), \dots, \bm{f}(\bm{x}_N) \right)^\top  \in \mathbb{R}^{N\times q}\\
\bm{R} & = R(\bm{x}, \bm{\tilde{x}})_{\bm{x},\bm{\tilde{x}}\in\bm{X}} \in \mathbb{R}^{N\times N}.
\end{split}
\end{equation}

The posterior mean $\mu$ and variance $\sigma^2$ can be thought of as the surrogate model's prediction and the corresponding uncertainty, and can be expressed as follows: 
\begin{equation} 
\label{eq:gp_prediction}
\begin{split}
{\mu}(\bm{x}) & = \bm{f}(\bm x)^\top \bm{\beta} + \bm{r}(\bm{x})^\top \bm{R}^{-1} (\bm{Y} - \bm{F}\bm{\beta}) \\
{\sigma}^2(\bm{x}) & = \sigma^2_z \left( 1 - \bm{r}(\bm{x})^\top \bm{R}^{-1} \bm{r}(\bm{x}) \right).
\end{split}
\end{equation}

\subsection{Multi-fidelity Gaussian processes}
The accuracy of numerical simulation outputs is often achieved at the expense of increased computational costs. Nevertheless, faster approximations using a coarser mesh, less complex models, or entirely different physical models can provide useful information, such as overall trends in the domain. Kennedy and O’Hagan~\cite{kennedy_predicting_2000} proposed the formulation presented here:
\begin{equation}
\label{eq:mf_kennedy}
y_{HF}(\bm{x}) = \rho y_{LF}(\bm{x}) + \delta(\bm{x}) \quad \text{s.t.} \quad y_{LF} \perp \delta.
\end{equation}
The high-fidelity output is expressed as the low-fidelity output scaled by a factor $\rho$, plus a discrepancy function $\delta$. The discrepancy function is modeled by a Gaussian process with the task of capturing the difference between 2 sequential levels of fidelity. Following this approach,~\citet{le_gratiet_multi-fidelity_2013} proposed a formulation where the LF surrogate model $\hat{y}_{LF}$ is used as a basis function to construct the HF surrogate model $\hat{y}_{HF} $ as follows:
\begin{equation}
\label{eq:mf_legratiet}
\hat{y}_{HF}(\bm{x}) = \sum_{k=1}^{q}{\beta_k f_k(\bm{x})} + \beta_{\rho} \hat{y}_{LF}(\bm{x}) + \delta(\bm x).
\end{equation}
This formulation can be extended to $L$ levels. This model is built recursively, starting from the lowest level to the highest. Each level $\ell$ is defined by their respective process variance $\sigma_{z,\ell}$, correlation length $\bm{\theta}_\ell$, and training data $\bm{D}_\ell$. Instead of using a single kernel that aggregates the spatial correlation of all levels, one kernel $R_\ell \in \mathbb{R}^{N_\ell \times N_\ell}$ is constructed for each level. Here, $N_\ell$ represents the number of samples of a given fidelity $\ell$. This formulation requires the use of a nested DoE : the sample locations of a given level must be present in all lower DoE ($\bm{D}_\ell \subseteq \bm{D}_{\ell-1}$). The mean and variance predictions at a query point $\bm x$ and fidelity level $\ell \geq 2$ are given respectively by:
\begin{align}
\mu_\ell(\bm{x}) & = \rho_{\ell-1}(\bm{x})\mu_{\ell-1}(\bm{x}) + \bm{f}_{\ell}(\bm{x})^\top \bm{\beta}_\ell + \bm{r}_\ell(\bm{x})^\top \bm{R}^{-1}_\ell \Big( \bm{Y}_\ell - \rho_{\ell-1}(\bm{D}_\ell) \odot \bm{Y}_{\ell-1}(\bm{D}_\ell) - \bm{f}_\ell(\bm{D}_\ell)^\top \bm{\beta}_\ell \Big)\\
\sigma^2_\ell(\bm{x}) & = \rho^2_{\ell-1}(\bm{x})\sigma^2_{\ell-1}(\bm{x}) + \underbrace{\sigma^2_{z,\ell} \left( 1 - \bm{r}_\ell(\bm{x})^\top \bm{R}^{-1}_\ell \bm{r}_\ell(\bm{x}) \right)}_{\sigma^2_{\delta, \ell}}.
\end{align}
In this paper, the regression function is constant $\bm{f}_\ell(\bm{x}) \to f_\ell = 1$, the scaling factor is a scalar value $\rho_\ell(\bm{x}) \to \rho_\ell$ and $r_\ell(\bm{x})$ is the correlation matrix between the training points $\bm{D}_\ell$ and the query point $\bm{x}$ : $\bm{r}_\ell(\bm{x}) = R(\tilde{\bm{x}}, \bm{x})_{\tilde{\bm{x}}\in \bm{D}_\ell}$. $\rho_{\ell-1}(\bm{D}_\ell)$ and $\bm{Y}_{\ell-1}(\bm{D}_\ell)$ denote the scaling factors and function values of the previous level for points in $\bm D_{\ell-1}$ and $\bm D_{\ell}$. The previous level function values are present in the previous and current level DOE. The mean and variance prediction functions for $\ell = 1$ reduce to Eq.~\eqref{eq:gp_prediction}. The model parameters $\{\bm{\beta}_\ell, \sigma_{z,\ell}, \bm{\theta}_\ell\}^L_{\ell=1}$, as well as the scaling factors $\{\rho_\ell\}^{L-1}_{\ell=1}$, are estimated by maximizing the restricted log-likelihood sequentially, starting from the lowest level to the highest. Once the hyper-parameters of a given level are estimated, they are held fixed for the subsequent levels. The variance contribution of a given level $\ell$ scaled up to the highest level $L$ can be computed using: 
\begin{equation}
\label{eq:variance_contribution}
\sigma^2_{\text{cont}}(\ell, \bm{x}) = \sigma^2_{\delta, \ell}(\bm{x}) \prod^{L-1}_{j=\ell}\rho^2_j
\end{equation}
This multi-fidelity GP formulation is used to construct surrogate models for the objective and the constraint functions employed in the Bayesian optimization approach discussed in Sections~\ref{sec:3} and~\ref{sec:4}.

\section{Multi-fidelity Bayesian optimization}\label{sec:3}
\subsection{Bayesian optimization}
Constrained Bayesian optimization employs Gaussian processes to model both the objective and constraint functions. At the beginning of the optimization process, the initial models are constructed using a limited number of function samples. At a given iteration $i$, the next sampling location is selected by optimizing an acquisition function $\alpha$:
\begin{equation}
\label{eq:constrained_acq_func}
\begin{aligned}
\bm{x}_{i} \in  \mathop{\arg \max}_{\bm{x} ~\in ~\mathbb{R}^d} & \quad \alpha(\bm{x}) \\
\text{s.t.} \quad
& \bm \mu_{g}(\bm{x}) \leq \bm 0\\
& \bm \mu_{h}(\bm{x}) = \bm 0\\
\end{aligned}
\end{equation}
where $\bm\mu_g(\bm{x})$ and $\bm\mu_h(\bm{x})$ are the mean predictions of the inequality and equality constraints. The SEGO framework~\cite{SasenaExplorationmetamodelingsampling2002} utilizes the Expected Improvement (EI) function~\cite{marchuk_bayesian_1975} as an acquisition function denoted by: 
\begin{equation}
\label{eq:ei}
\alpha(x)=\text{EI}(\bm{x}) = 
\begin{cases}
(f_{\min} - \mu_f(\bm{x}))\Phi\left(\frac{f_{\min}-\mu_f(\bm{x})}{\sigma_f({\bm{x})}}\right) + \sigma_f(\bm{x})\phi\left(\frac{f_{\min}-\mu_f(\bm{x})}{\sigma_f(\bm{x})}\right) &\text{if } \sigma_f(\bm{x}) > 0 \\
0 &\text{if } \sigma_f(\bm{x}) = 0
\end{cases}
\end{equation}
The output mean $\mu_f(\bm{x})$ and variance $\sigma_f^2(\bm{x})$ of the objective GP are used to identify locations with probabilities of improving the current best feasible minimum objective sampled value $f_{\min}$. If no point is feasible, the objective value with the least constraint violation is used instead. $\Phi$ and $\phi$ are, respectively, the cumulative distribution and the probability density functions of \(\mathcal{N}(0, 1)\). The blackbox functions are sampled at the infill location $\bm{x}_{i}$ and the corresponding objective and constraint data are added to the DoE\mycomment{ : $\bm{\text{DoE}} = \{\bm{x}_l, f(\bm{x}_l), \bm{g}(\bm{x}_l), \bm{h}(\bm{x}_l)\}_{l=1,\dots,N}$}. Even if the sampled point is unfeasible, it is still added to the GP training points to improve the accuracy of the objective and constraints GPs. This process is repeated until a convergence criterion is met, such as a maximum number of iterations. The SEGO framework is summarized by Algorithm~\ref{alg:sego_algo}. \mycomment{$\Omega_c := \{\bm{x}\in \mathbb{R}^d \; | \; \bm{g}(\bm{x}) \leq \bm{0}, \; \bm{h}(\bm{x}) = \bm{0}\}$ is used to denote the feasible domain.}

\begin{algorithm}
\caption{The Super Efficient Global Optimization (SEGO) framework.}
\label{alg:sego_algo}
\begin{algorithmic}
\State \textbf{Input:}~ Design space $\Omega$, objective and constraints oracles $ \{f, \bm{g}, \bm{h}\}$, and budget, i.e., $ \text{max\_iter}$.
\State \textbf{Output:}~ $f_{\min}$ the best feasible point in terms of the objective function $f$.

\State Generate an initial DoE; 

\State Set $ f_{\min}$ to the best feasible point in terms of the objective function $f$;
\State $i \gets 0$;
\While{$ i \le \text{max\_iter} $}
    \State Build the surrogate models using GPs for the objective and constraint functions;
    \State $ \bm{x}_{i} \in \displaystyle \mathop{\arg\max}_{\bm{x}\in \Omega } \{ \alpha (\bm{x}) \quad \text{s.t.} \quad \bm\mu_g(\bm{x}) \leq \bm{0} , \; \bm\mu_h(\bm{x}) = \bm{0} \}$; \Comment Find infill location
    \State Evaluate $f$, $\bm{g}$ and $\bm{h}$ at $\bm{x}_i$;
    \State Update the DoE;
 \State Set $ f_{\min}$ to the best feasible point in terms of the objective function $f$;
        \State $i \gets i+1$;
\EndWhile

\end{algorithmic}
\end{algorithm}

\subsection{Bayesian optimization with only objective multi-fidelity selection criteria}
Multi-fidelity Bayesian optimization leverages multiple levels of fidelity to reduce the cost of optimization and increase the precision of the surrogate models with computationally cheaper data. The MFSEGO framework employed in~\cite{meliani_multi-fidelity_2019, charayron2023towards} is heavily derived from the SEGO framework previously discussed. The objective and constraint functions are modeled using the Le Gratiet Multi-Fidelity Kriging model, thus requiring a nested DoE such that $\mbox{DoE}_{\ell} \subseteq \mbox{DoE}_{\ell-1}$ for $\ell=2,\dots,L$.  The infill location and the fidelity level at which to sample said location are computed sequentially as a two-step approach. For the first step, the infill location is selected using an acquisition function such as EI~\cite{marchuk_bayesian_1975} or the scaled WB2~\cite{bartoli_adaptive_2019}. $f_{\min}$ corresponds to the best feasible objective evaluated at the highest level of fidelity. From Eq.~\eqref{eq:variance_contribution}, the reduction in variance of the surrogate model, if the function were to be sampled at level $\ell$ for the objective function $f$, can be computed as follows: 
\begin{equation}
\label{eq:variance_reduction}
\sigma_{(\text{red},0)}^2(\ell, \bm{x}_{i}) = \sum_{\ell'=1}^{\ell} \sigma_{(\text{cont},0)}^2(\ell', \bm{x}_{i})= \sum_{\ell'=1}^{\ell}\sigma^2_{(0, \delta, \ell')}(\bm{x}_{i}) \prod_{j=\ell'}^{L-1}{\rho^2_{(j,0)}}.
\end{equation}
The subscript $(\mbox{red},0)$ indicates that the reduction is related to the objective function.
The nested DoE requires that if the function is sampled at fidelity $\ell$, it must also be sampled at all lower fidelity levels. The total cost of sampling a level is given by $\sum_{\ell'=1}^{\ell}{c_{\ell'}}$
where $c_{\ell'}$ represents the individual cost of each level $\ell'$ for evaluating the objective and constraint functions. 
To select the fidelity level for the second step, we use the best ratio of variance reduction to the cost squared:
\begin{equation}
\label{eq:obj_only}
\ell^{\mbox{obj}}_{i} = \mathop{\arg\max}_{\ell \in \{1, \dots, L
\}} \frac{\sigma^2_{(\text{red},0)}(\ell, \bm{x}_{i})}{\left(\sum_{\ell'=1}^{\ell}{c_{\ell'}}\right)^2}.
\end{equation}
\citet{meliani_multi-fidelity_2019} suggest that the cost should be squared, since the variances are scaled with the squared scaling factor.

\section{MFSEGO using both objective and constraints multi-fidelity criteria}
\label{sec:4}

Previous work only applied the fidelity criterion to the objective surrogate model, e.g., \cite{meliani_multi-fidelity_2019,charayron2023towards}. Since constraints can also be modeled using different fidelity levels, in this section, three other fidelity criteria are proposed to incorporate the variances of the constraint models into the fidelity level selection process. The fidelity criterion previously discussed (see Eq.~\eqref{eq:obj_only}) will be referred to as the ``objective-only'' fidelity criterion.

Inspired by Eq.~\eqref{eq:variance_reduction}, for a given \( k \in \{0, \ldots, m+p\} \), the reduction in variance of the \( k \)-th GP component is given by  
\begin{equation}
\label{eq:variance_reduction:k}
\sigma_{(\text{red},k)}^2(\ell, \bm{x}_{i}) = \sum_{\ell'=1}^{\ell} \sigma_{(\text{cont},k)}^2(\ell', \bm{x}_{i})= \sum_{\ell'=1}^{\ell}\sigma^2_{(k, \delta, \ell')}(\bm{x}_{i}) \prod_{j=\ell'}^{L-1}{\rho^2_{(j,k)}}.
\end{equation}
Let $\sigma_{(\text{norm},k)}^2$ be the normalized reduction in variance of \( g_k \) with respect to the cost:
\begin{equation}
\label{eq:variance_reduction:norm}
\sigma_{(\text{norm},k)}^2(\ell, \bm{x}_{i}) = \sigma_{(\text{red},k)}^2(\ell, \bm{x}_{i}) \left(\displaystyle \sum_{\ell'=1}^{\ell}{c_{\ell'}}\right)^{-2}.
\end{equation}
The ``Average'' variance reduction criterion averages the variance reduction of the objective and constraint models for a level $k$ and divides it by the level's total cost squared. The level with the best ratio is selected. This approach can be extended to $m+p+1$ GPs (one GP model for the objective function and $m+p$ GPs for the constraints) as written: 
\begin{equation}
\label{eq:average}
\ell^{\mbox{avg}}_i = \mathop{\arg\max}_{\ell \in \{1, \dots, L\}} \displaystyle \sum_{k=0}^{m+p}   \sigma_{(\text{norm},k)}^2(\ell, \bm{x}_{i}).
\end{equation}
The ``Optimistic'' variant computes the best level $\ell^{\mbox{min}}_i$ for each GP and selects the overall lowest from the $m+p+1$ GPs: 
\begin{equation}
\label{eq:optimistic}
\ell^{\mbox{min}}_i  =  \min_{k \in \{0, \dots, m+p\}} \left(\mathop{\arg\max}_{\ell \in  \{1, \dots, L\}}
\sigma_{(\text{norm},k)}^2(\ell, \bm{x}_{i}) \right).
\end{equation}
The third variant, referred to as the ``Pessimistic'' one,  computes the overall highest level from the $m+p+1$ GPs: \mycomment{related to the constraints and the objective functions.}

\begin{equation}
\label{eq:pessimistic}
\ell^{\mbox{max}}_i  = \max_{k \in \{0, \dots, m+p\}} \left(\mathop{\arg\max}_{\ell \in  \{1, \dots, L\}} \sigma_{(\text{norm},k)}^2(\ell, \bm{x}_{i}) \right).
\end{equation}
The MFSEGO framework using both objective and constraints multi-fidelity criteria is given by Algorithm~\ref{alg:mfsego}, where $f$ represents the objective function, $\bm{g}$ the inequality constraint functions, and $\bm{h}$ the equality constraint functions.

\begin{algorithm}
\caption{The multi-fidelity based SEGO (MFSEGO) for constrained optimization.}\label{alg:mfsego}
\begin{algorithmic}
\State \textbf{Input:}~ Design space $\Omega$, objective and constraints oracles per fidelity level $ \{f, \bm{g}, \bm{h}\}_{\ell\in\{1,\dots,L\}}$, cost per fidelity level $\{c_\ell\}_{(\ell\in\{1,\dots,L\})}$, and a budget, i.e., $ \text{max\_iter}$.
\State \textbf{Output:}~ $f_{\min}$ the best feasible point in terms of the high-fidelity objective value.
\State Generate an initial nested $ \{\mbox{DoE}_\ell\}_{\ell \in \{1,\ldots,L\}}$;
\State Set $ f_{\min}$ to the best feasible point in terms of the objective function in $\mbox{DoE}_{\ell=L}$;
\State $i \gets 0$;
\While{$ i \le \text{max\_iter} $}
    \State Build  GPs for the objective and the constraints;
    \State $ \bm{x}_{i} \in \displaystyle \mathop{\arg\max}_{\bm{x} \in \Omega} \{ \alpha (\bm{x}) \quad \text{s.t.} \quad \bm\mu_g(\bm{x}) \leq \bm{0} , \; \bm\mu_h(\bm{x}) = \bm{0} \}$; \Comment Find infill location

    \State Select a fidelity  level $\ell_i \in  \{\ell^{\mbox{obj}}_i, \ell^{\mbox{avg}}_i, \ell^{\mbox{min}}_i, \ell^{\mbox{max}}_i \}$;

    \For {$ \ell' \; \text{in} \; \{1, \dots, \ell_i \}$}    \Comment Sample and update each fidelity level up to $ \ell_i $
        \State Evaluate the objective and constraints at $\bm{x}_i$;
        \State Update the nested $\mbox{DoE}_{\ell'}$;
    \EndFor

\State Update $ f_{\min}$ as the best feasible point in terms of the objective value in $\mbox{DoE}_{\ell=L}$;
   \State $i \gets i+1$
    
\EndWhile


\end{algorithmic}
\end{algorithm}

\section{Numerical results}\label{sec:5}
\subsection{Implementations details}
All surrogate models are built using the MFK model implemented in the SMT package~\cite{saves_smt_2024}. The GP hyperparameters are estimated by maximizing the restricted log-likelihood. The optimized hyperparameters from the previous iteration are used as a starting location for the sub-optimization problem in addition to $3d$ other starting locations generated using LHS.
%
The constrained acquisition function is optimized using SLSQP~\cite{Kraft1988} because of its ability to handle both inequality and equality constraints. A multi-start strategy is used to escape any local minima. The starting locations are generated using LHS. As the number of infill samples increases, the output variance of the objective surrogate model tends to decrease. The gradient of $\text{EI}$ can quickly fall under the machine epsilon $\epsilon$, making it practically difficult to optimize numerically. To address this issue, Ament et al.~\cite{ament2023unexpected} proposed a numerically stable implementation of the Log Expected Improvement acquisition function as denoted by: 

\begin{equation}
\label{eq:log_ei}
\alpha(\bm{x})=\log\text{EI}(\bm{x}) = \mathtt{log\_h}\left(\frac{f_{\min} - \mu_f(\bm{x})}{\sigma_f(\bm{x})}\right) + \log(\sigma_f(\bm{x}))
\end{equation}
where
 \[
\mathtt{log\_h}(z) = 
\begin{cases}
\log(\phi(z) + z\Phi(z)) & z > -1 \\
-\frac{z^2}{2} - c_1 + \mathtt{log1mexp}(\log(\mathtt{erfcx}(-z/\sqrt{2})|z|) + c_2) & -1/\sqrt{\epsilon} < z \leq -1 \\
-\frac{z^2}{2} - c_1 - 2\log(|z|) & z \leq -1/\sqrt{\epsilon}
\end{cases}
\]
\noindent
with $c_1 = \frac{1}{2}\log(2\pi)$, $c_2 = \frac{1}{2}\log(\pi/2)$ and $\epsilon>0$ is set to the epsilon machine tolerance. The functions \texttt{log1mexp} and \texttt{erfcx} are defined as $\texttt{log1mexp}(z) = \log(1 - \exp(z))$ and $\texttt{erfcx}(z) = \exp(z^2) \,\text{erfc}(z)$. 
This $\log\text{EI}$ function returns non-numerically zero values even when the probability of improvement would lead traditional $\text{EI}$ to fall below epsilon.\mycomment{decreases considerably.}
Furthermore, when the output variances of the GPs decrease, so does the variance reduction $\sigma^2_{\text{red}}$. To avoid $\sigma^2_{(\text{norm},k)}$ falling below machine epsilon, the cost of every level $c_{\ell=\{1,\dots,L\}}$ is normalized by the highest level $c_{\ell=L}$.

\subsection{Results on analytical test problems}

The MFSEGO algorithm, as described in Sections~\ref{sec:2} and~\ref{sec:3}, is applied to multi-fidelity constrained problems: the Rosenbrock problem (noted MF Rosenbrock)~\cite{lam_multifidelity_2015}, the Branin problem (noted MF Branin)~\cite{qian_optimization_2021}, the Sasena problem (noted MF Sasena)~\cite{qian_optimization_2021}, and the Gano problem (noted MF Gano)~\cite{Gano2005}.
Table~\ref{tab:problem_properties} summarizes the test problem characteristics. The respective objective and constraint functions are given in Appendix~\ref{sec:appendixB}. The objective and constraint functions are each defined with two fidelity levels. The LF functions are obtained by adding a perturbation to their HF counterpart. The global minimum of each problem is located on the constraint's boundary. A set of 25 initial feasible DoE was created for each problem.

\begin{table}[ht]
\centering
\caption{Analytical test problems properties.}
\label{tab:problem_properties}
\begin{tabular}[t]{lcccccccc}
\toprule
Problem & Dimensions & Domain & $ \min\limits_{\bm{x} \in \Omega_g} \; f $ & $ \arg \! \min\limits_{\bm{x} \in \Omega_g} \; f $ & $ \rho_{f} $ & $ \rho_{g} $ & LF DoE & HF DoE\\
\midrule
MF Rosenbrock & 2 & $[-2, \; 2]$ & 0.1785 & $[0.5777, \; 0.3325]$ & $\approx1$ & $0.9986$ & 6 & 3 \\
MF Branin & 2 & $[0, \; 1]$ & 5.5757 & $[0.9676, \; 0.2067]$ & $\approx1$ & $0.8163$ & 6 & 3 \\
MF Sasena & 2 & $[0, \; 5]$ & -1.1743 & $[2.7450, \; 2.3523]$ & $0.3513$ & $0.2981$ & 6 & 3 \\
MF Gano & 2 & $[0.1, 10]$ & 5.6684 & $[0.8842, 1.1507]$ & $\approx1$ & $0.9723$ & 6 & 3 \\
\bottomrule
\end{tabular}
\end{table}

\begin{table}[ht]
\centering
\caption{Multi-fidelity framework comparison}
\begin{tabular}[t]{lcccc}
\toprule
Framework & SMT Surrogate & DoE & MF constraint & Equality constraints \\
\midrule
MFSEGO & MFK & Nested & Sub-solver & \checkmark \\
VF-EI & MFCK & Non-nested & PoF & x \\
VF-PI & MFCK & Non-nested & PoF & x \\
\bottomrule
\label{tab:mf_framework_comparison}
\end{tabular}
\end{table}

\begin{figure}[h]
\includegraphics[scale=1]{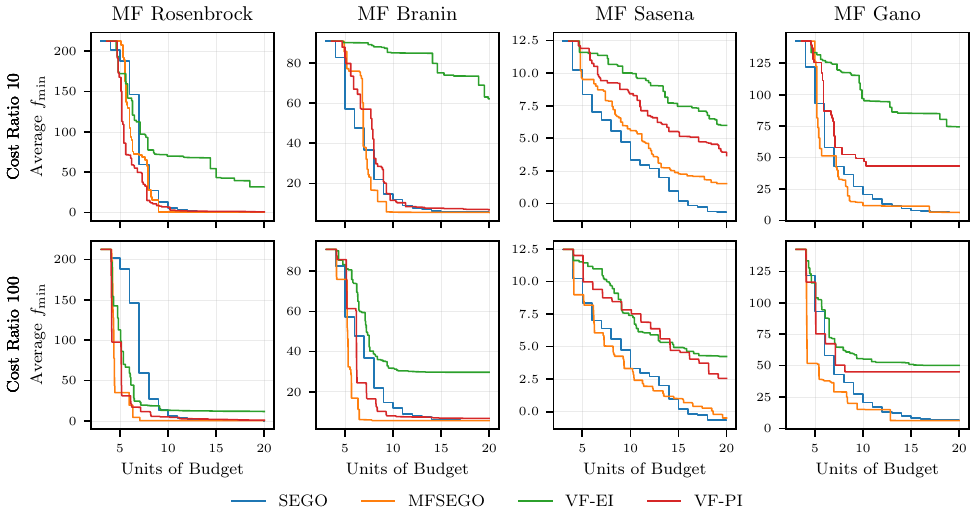}
\caption{Multi-fidelity BO frameworks comparison. Analytical problems convergence results for the Rosenbrock, Branin, Sasena and Gano problems with two different cost ratios. Each optimization framework was run on 25 DoE.}
\label{fig:mf_frameworks_comparison}
\end{figure}

Figure~\ref{fig:mf_frameworks_comparison} compares MFSEGO with SEGO and two other multi-fidelity constrained BO frameworks: VF-EI~\cite{Zhang2018} and VF-PI~\cite{Ruan2020}. These two approaches do not require the use of a nested DoE. To handle the constraints when maximizing their respective acquisition function, both frameworks penalize it with the probability-of-feasibility (PoF)~\cite{Schonlau1997}. This method is not well suited for handling multiple constraints as it often drives the acquisition function to be numerically zero over the feasible domain and cannot efficiently handle equality constraints. Furthermore, both frameworks were implemented with the MFCK model in the SMT package, as it supports non-nested DoE. Their acquisition function is maximized using L-BFGS-B~\cite{Byrd1995} with a multi-start strategy. Table~\ref{tab:mf_framework_comparison} summarizes the properties of MFSEGO, VF-EI, and VF-PI. All methodologies are available in the BOMA repository~\footnote{\hyperlink{https://github.com/oihanc/boma}{https://github.com/oihanc/boma}}. Figure~\ref{fig:mf_frameworks_data_profile} in Appendix~\ref{sec:appendix_data_profile} displays the data profile for the same experiments. It illustrates that MFSEGO converges more efficiently than VF-EI and VF-PI. VF-PI is shown to be more efficient than SEGO except for a few instances where it failed to converge.

\begin{figure}[h]
\includegraphics[scale=1]{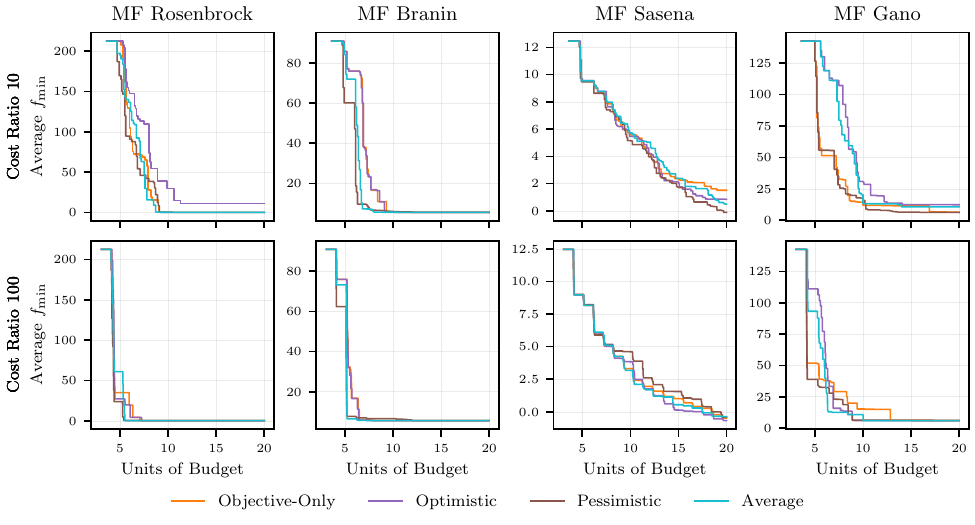}
\caption{Fidelity criteria comparison. Analytical problems convergence results for the MF Rosenbrock, MF Branin, MF Sasena and MF Gano problems with two different cost ratios. Each fidelity criterion was run on 25 DoE.}
\label{fig:test_problems_convergence_fid_crit}
\end{figure}

Figure~\ref{fig:test_problems_convergence_fid_crit} compares the fidelity criteria as discussed in Section~\ref{sec:3}. When applied to the 4 previously mentioned test problems, no fidelity criterion distinguished itself from the others. It should be noted that the Objective-Only criterion corresponds to the state of the art MFSEGO algorithm. Moreover, the difference in the required budget for convergence between the proposed fidelity criteria decreases when using a higher cost ratio.

\subsection{Multi-fidelity aero-structural wing design optimization}

The fidelity criteria proposed in Section~\ref{sec:3} were benchmarked on a aero-structural test case using OpenAeroStruct~\cite{Jasa2018}. Table~\ref{tab:oas_definition} summarizes the problem characteristics. The objective is to minimize the fuel burn using the Breguet equation: 
\begin{equation}
\mathtt{fuel burn} = (W_0 + W_s) \left( \exp(\frac{R C_T}{a M} \frac{C_D}{C_L}) - 1 \right)
\end{equation}
with regards to the angle-of-attack, the wing twist, and the thickness of the spar. $W_0$ and $W_s$ respectively symbolize the wing empty weight without fuel and structural elements, and the wing structural weight due to the spar. $R$ corresponds to the aircraft total range; $C_T$ to the specific fuel consumption; $a$ to the speed of sound and $M$ to the Mach number. Both the wing twist and the spar thickness are defined using 5 control points uniformly distributed along the wing, from which a B-Spline is fitted. Two constraints are imposed: the lift force must be equal to the weight, and the wing structural loads must respect a safety factor of $2.5$. It should be noted that the fuel burn objective and the lift-to-weight constraint are well correlated, contrary to the structural failure constraint.

\begin{table}[h!]
    \centering
    \caption{Definition of the wing aerostructural optimization problem.}
    \label{tab:oas_definition}
    \begin{tabular}{@{}lccccc@{}}
        \toprule
        & \textbf{Function/variable} & \textbf{Description} & \textbf{Quantity} & \textbf{Range} & $\rho_\text{Pearson}$ \\
        \midrule
        \textbf{Minimize} & $\mathtt{fuel\_burn}$      & Fuel Burn  & 1  & -- & $0.99$   \\
        \midrule
        \multirow{1}{*}{with respect to} 
                                    & AoA & Angle-of-Attack [deg] & 1 & [8, 12] & -- \\
                                    & $\alpha_{\text{twist}}$ & Twist [deg] & 5 & $[-6, 3]$ & -- \\
                                    & $t_\text{spar}$ & Spar thickness [m] & 5 & [0.0015, 0.05] & -- \\
        \midrule
        \text{subject to}           &   $L/W = 1$    &   Lift-to-Weight ratio & 1 & -- & $0.98$\\
                                    & $\sigma_\text{VM} \leq 0$ & Von-Mises criterion & 1 & -- & $0.70$ \\
        \midrule
        \bottomrule
    \end{tabular}
\end{table}

\begin{figure}[h]
\centering
\includegraphics[scale=0.8]{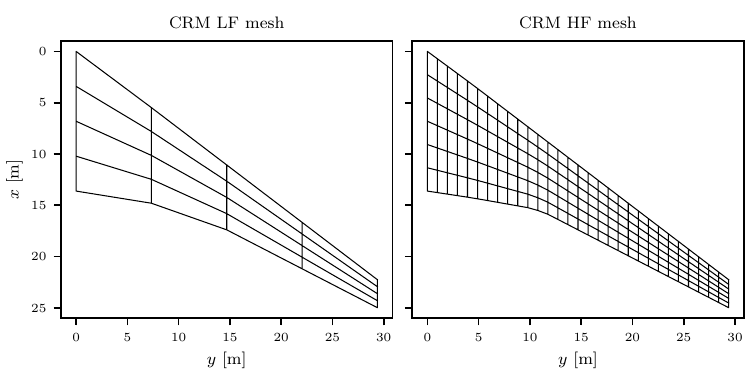}
\caption{CRM multi-fidelity VLM mesh comparison. The LF and HF meshes use $4 \times 4$ panels and $6 \times 30$ panels respectively. The cost ratio between the 2 levels is 30.}
\label{fig:oas_mf_meshes}
\end{figure}

Figure~\ref{fig:oas_mf_meshes} illustrates the low- and high-fidelity meshes used for this test case. The aerodynamic forces are computed using the Vortex Lattice method (VLM) which assumes an incompressible and inviscid fluid flow. The structural forces are computed with a finite-element method. The 4 fidelity criteria were benchmarked on 25 instances, each starting with a different nested DoE comprised of 122 LF points and 12 HF points. Due to the equality constraint, no initial DoE was feasible. The cost ratio between these fidelities is 30. Using this value led the optimization runs to rarely sample the HF level. Excessive LF sampling increases the BO overhead cost; thus, to encourage more frequent HF evaluations, the cost ratio used in Eq.~\eqref{eq:variance_reduction:norm} was reduced to 10.

\begin{figure}[h]
\includegraphics[scale=1.0]{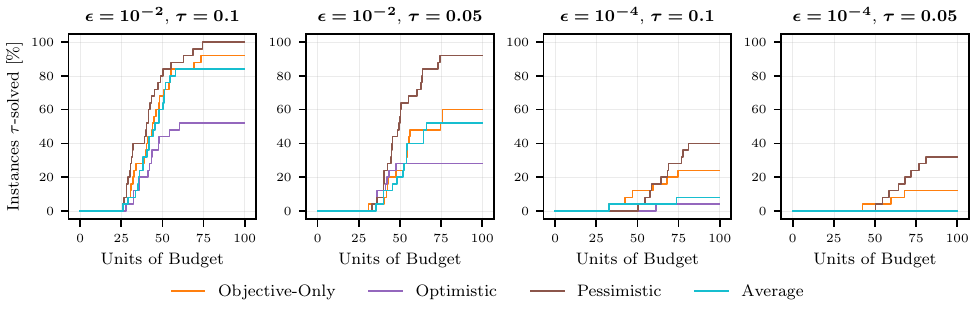}
\caption{Fidelity criteria comparison on the wing aerostructural optimization problem. $\epsilon$ and $\tau$ denote the absolute constraint tolerance and the relative objective tolerance respectively. Each fidelity criterion was run on 25 DoE.}
\label{fig:oas_data_profiles}
\end{figure}

Figure~\ref{fig:oas_data_profiles} illustrates data profiles for various absolute constraint tolerance $\epsilon$ and relative objective tolerance $\tau$. Here $\epsilon$ is defined as the Root Square Constraint Violation (RSCV) defined by:
\begin{equation}
\text{RSCV} = \sqrt{\sum^{m}_{i=1} \max(g_i, 0)^2 + \sum^p_{j=1} h_j^2}.
\end{equation}
Additionally, a run is said to be $\tau$-solved if the following criterion is met:
\begin{equation}
f^i - f_0 \ge (1-\tau)(f^\star - f_0)
\end{equation}
where $f_0$ is the greatest first feasible objective value over all the runs, and $f^\star$ is the best feasible objective value over all the runs.

\begin{figure}[h]
\centering
\includegraphics[scale=1.0]{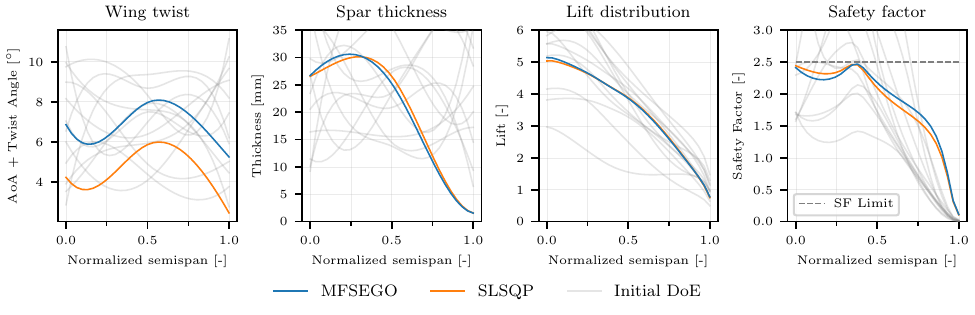}
\caption{Comparison of the best solutions found with MFSEGO and SLSQP for wing twist, tubular spar thickness, lift distribution and structural safety factor.}
\label{fig:oas_twist_thickness_lift_failure}
\end{figure}

The Pessimistic criterion clearly outperforms the other fidelity criteria, especially when imposing stricter constraints and objective tolerances. Compared to the state of the art approach using the Objective-Only criterion, the Pessimistic fidelity criterion solved 9\% to 200\% more instances, with the difference growing for stricter tolerances. Figure~\ref{fig:oas_twist_thickness_lift_failure} compares the wing twist, the spar thickness, the lift distribution, and the structural safety factor along the wing semispan of the best solution obtained using MFSEGO, and  SLSPQ. The corresponding configurations from the  initial HF DOE are also shown. As a gradient-based solver, SLSQP requires a single starting point, which was selected as the point with the smallest constraint violation in the HF DOE. Being a mono-fidelity solver, SLSQP operated exclusively on the HF objective and constraint functions. An optimal lift distribution would follow an elliptical shape. In this case, both MFSEGO and SLSQP solutions are limited by the structural constraint. 

\section{Conclusions and perspectives}\label{sec:6}

In this paper, we investigate the potential of incorporating multiple fidelity levels for both the constraints and the objective function into a multi-fidelity constrained Bayesian optimization framework, i.e., MFSEGO. Different criteria are considered to handle the multi-fidelity aspects, ranging from the most optimistic approach, where decisions are mainly based on the best affordable fidelity level, to the most pessimistic approach, where decisions are guided by the worst affordable fidelity level. 

MFSEGO, VF-EI~\cite{Zhang2018}, and VF-PI~\cite{Ruan2020} were benchmarked on 4 constrained multi-fidelity test problems. SEGO was also benchmarked as a mono-fidelity reference. Overall, MFSEGO proved to be more efficient. Three proposed multi-fidelity criteria were benchmarked against the existing objective-only multi-fidelity criterion and showed comparable performance on the analytical test problems. However, the pessimistic criterion proved more efficient on the aero-structural wing design problem: for a limited budget and strict tolerances, the pessimistic criterion solved more instances than any other multi-fidelity criterion. Results indicate that multi-fidelity approaches can reduce the computational budget while achieving similar convergence. Future work includes integrating mixed-integer variables and fidelity-dependent design spaces to evaluate the full potential of the MFSEGO approach on a more complex benchmark~\cite{BaDiMoLeSa2023}, such as a test case similar to the Bombardier Research Aircraft Configuration~(BRAC)~\cite{priem_2020}.

\section{Appendix}\label{sec:appendix}


\subsection{Rosenbrock analytical test problem}\label{sec:appendixB}
The Rosenbrock multi-fidelity test problem~\cite{lam_multifidelity_2015} is defined by:
\begin{equation}
    \label{eq:rosenbrock}
    \begin{split}
        f_{HF} & = (1 - x_0)^2 + 100(x_1 - x_0^2)^2\\
        f_{LF} & = f_{HF}(x_0, x_1) + 0.1\sin(10 x_0 + 5 x_1)\\
        g_{HF} & = x_0^2 + \sqrt{x_1 - 1} \leq 0 \\
        g_{LF} & = g_{HF}(x_0, x_1) - 0.1\sin(10 x_0 + 5 x_1) \leq 0.
    \end{split}
\end{equation}

\subsection{Branin analytical test problem}\label{sec:appendixC}
The modified Branin multi-fidelity test problem~\cite{qian_optimization_2021} is defined by:
\begin{equation}
    \label{eq:branin}
    \begin{split}
        f_{HF} & = \left( 15 x_1 - \frac{5.1}{4 \pi^2} (15 x_0 - 5)^2 + \frac{5}{\pi} (15 x_0 - 5) - 6 \right)^2 \\ 
        &\quad + 10 \left( 1 - \frac{1}{8 \pi} \right) \cos(15 x_0 - 5) + 10 + 5 x_0 \\
        f_{LF} & = f_{HF}(x_0, x_1) - \cos(0.5x_0) - x_1^3\\
        g_{HF} & = -x_0 x_1 + \frac{1}{5} \leq 0 \\
        g_{LF} & = -x_0 x_1 +\frac{3}{10} x_0 - \frac{7}{10}x_1 \leq 0.
    \end{split}
\end{equation}

\subsection{Sasena analytical test problem}\label{sec:appendixD}
The Sasena multi-fidelity test problem~\cite{qian_optimization_2021} is defined by:
\begin{equation}
    \label{eq:sasena}
    \begin{split}
        f_{HF} & = 2 + 0.01(x_1 - x_0^2)^2 + (1 - x_0)^2 + 2(2 - x_1)^2 + 7\sin( x_0)\sin(0.7 x_0 x_1) \\
        f_{LF} & = f_{HF}(x_0, x_1) + \exp{x_0} - x_1^3 \\
        g_{HF} & = -\sin(x_0 - x_1 - \pi/8) \leq 0 \\
        g_{LF} & = g_{HF}(x_0, x_1) + 0.2 x_1 - 0.7 x_0 + x_0 x_1 \leq 0.
    \end{split}
\end{equation}

\subsection{Gano analytical test problem}\label{sec:appendixE}
The Gano multi-fidelity test problem~\cite{Gano2005} is defined by:
\begin{equation}
    \label{eq:gano}
    \begin{split}
        f_{HF} & = 4 x_1^2 + x_2^3 + x_1 x_2 \\
        f_{LF} & = 4 (x_1 + 0.1)^2 + (x_2 - 0.1)^3 + x_1 x_2 + 0.1 \\
        g_{HF} & = 1/x_1 + 1/x_2 - 2 \leq 0 \\
        g_{LF} & = 1/x_1 + 1/(x_2 + 0.1) - 2 - 0.001\leq 0. \\
    \end{split}
\end{equation}

\subsection{Multi-fidelity BO frameworks data profile}
\label{sec:appendix_data_profile}
\begin{figure}[h]
\centering
\includegraphics[scale=1]{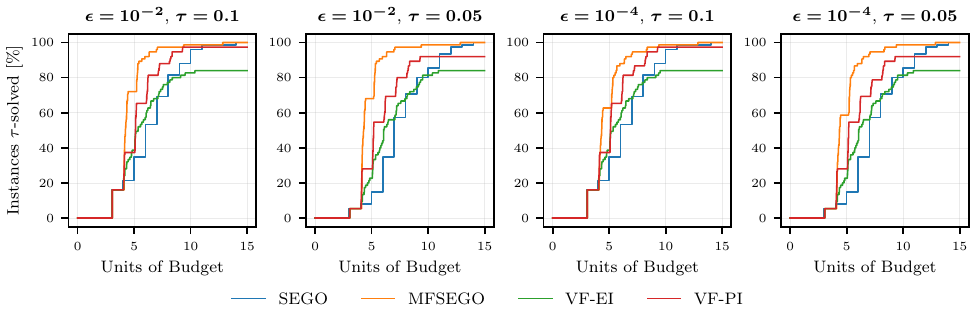}
\caption{Multi-fidelity BO frameworks comparison. Data profile generated with the MF Rosenbrock, MF Branin and MF Gano test problems with a fixed cost ratio of 100.}
\label{fig:mf_frameworks_data_profile}
\end{figure}
Figure~\ref{fig:mf_frameworks_data_profile} aggregates the optimization results displayed in Fig.~\ref{fig:mf_frameworks_comparison} for the MF Rosenbrock, MF Branin and the MF Gano test functions. The MF Sasena function was excluded due to its poor objective and constraint fidelity correlation factors. In decreasing order of convergence efficiency, the multi-fidelity framework ranking is: MFSEGO, VF-PI, and VF-EI. Given 20 units of budget, only MFSEGO and its mono-fidelity counterpart solved all 75 instances.




\bibliography{references}


\end{document}